# Frames and Bases in Tensor Product of Hilbert Spaces

Amir Khosravi

Faculty of Mathematical Sciences and Computer Engineering
University For Teacher Education, Taleghani Ave. 599
Tehran 15614, Iran
e-mail: khosravi@saba.tmu.ac.ir,   khosravi_amir@yahoo.com

M. S. Asgari

Department of Mathematics, Science and Research Branch
Islamic Azad University, Tehran, Iran
e-mail: msasgari@yahoo.com

## Abstract

In this article we develop a theory for frames in tensor product of Hilbert spaces. We show that like bases if $Y_1, \cdots, Y_n$ are frames for $H_1, \cdots, H_n$, respectively, then

$$Y_1 \otimes \cdots \otimes Y_n = \{y_1 \otimes \cdots \otimes y_n \ : \ y_1 \in Y_1, \cdots y_n \in Y_n\}$$

is a frame for $H_1 \otimes \cdots \otimes H_n$. Moreover we consider the canonical dual frame in tensor product space. We further obtain a relation between the dual frames in Hilbert spaces, and their tensor product.

## 1 Introduction

In 1946 Gabor [7] introduced a technique for signal processing which led eventually to wavelet theory. Later in 1952 Duffin and Schaeffer [5] introduced frame theory for Hilbert spaces. In 1986 Daubechies, Grossmann and Meyer [4] show that Duffin and Schaeffer's definition was an abstraction of Gabor's concept. Nowadays frames work as an alternative to orthonormal bases in Hilbert spaces which has many advantages [9]. Since tensor product is useful in approximation theory, in this article we consider the frames in tensor product of Hilbert spaces and extend some of the known results about bases to frames.

Let $H$ be a separable complex Hilbert space. As usual we denote the set of all bounded linear operators on $H$ by $B(H)$. We use $\mathbb{N}$, $\mathbb{Z}$, $\mathbb{R}$ and $\mathbb{C}$ to denote the natural numbers, integers, real numbers and complex numbers, respectively. $I$, $J$ and every $J_i$ will denote generic countable (or finite) index sets. We will always use $E_1 = \{e_i\}_{i \in I}$ and $E_2 = \{u_j\}_{j \in J}$ to denote orthonormal bases for $H$ and $K$, respectively. A sequence $\{x_n\}$ in a Hilbert space $H$ is called a **frame** for $H$, if there exist two constants $A$, $B > 0$ such that

$$A\|x\|^2 \leq \sum_n |<x, x_n>|^2 \leq B\|x\|^2, \qquad \text{for all } x \in H. \tag{1}$$

The numbers $A$ and $B$ are called the frame bounds. The frame $\{x_n\}$ is called a **tight frame** if we can choose $A = B$ and a normalized tight frame if $A = B = 1$. Therefore $\{x_n\}$ is a normalized tight frame if and only if for every $x \in H$, $\|x\|^2 = \sum_n |<x, x_n>|^2$. If $\{x_n\}$ is a normalized tight frame, then for every $x \in H$, $x = \sum_n <x, x_n> x_n$ (at least in the weakly convergent sense). Conversely, if $\{x_n\}$ is a sequence in $H$ such that the equation $x = \sum_n <x, x_n> x_n$ holds for every $x \in H$ (the convergence can be either in the weakly convergent sense or in the norm convergent sense) then $\{x_n\}$ is a normalized tight frame for $H$. Obviously every orthonormal basis is a normalized tight frame.

Throughout this paper, all of the Hilbert spaces will be separable and complex. For convenience we will denote the inner product of all Hilbert spaces by $<.,.>$.

## 2 Frame In Tensor Product

In this section we consider the tensor product of Hilbert spaces and we generalize some of the known results about bases to frames. There are several ways of defining the tensor product of Hilbert spaces. Folland in [6], Kadison and Ringrose in [11] have represented the tensor product of Hilbert spaces $H$ and $K$ as a certain linear space of operators. Since we used their results firstly we state some of the definitions.

Let $H$ and $K$ be Hilbert spaces. Then we consider the set of all bounded antilinear maps from $K$ to $H$. The operator norm of an antilinear map $T$ is defined as in the linear case:

$$\|T\| = \sup_{\|x\|=1} \|Tx\|. \tag{2}$$

The adjoint of a bounded antilinear map $T$ is defined by

$$<T^*x, y> = <Ty, x> \quad \text{for all } x \in H, y \in K. \tag{3}$$

Note that the map $T \longmapsto T^*$ is linear rather than antilinear. Suppose $T$ is an antilinear map from $K$ into $H$ and $E_1 = \{e_i\}_{i \in I}$ and $E_2 = \{u_j\}_{j \in J}$ are orthonormal bases for $H$ and $K$, respectively. Then by the Parseval identity

$$\sum_j \|Tu_j\|^2 = \sum_i \|T^*e_i\|^2 \tag{4}$$

This shows that $\sum_j \|Tu_j\|^2$ is independent of the choice of basis $E_2$.

**Definition 2.1** *Let $H$ and $K$ be Hilbert spaces. Then the tensor product of $H$ and $K$ is the set $H \otimes K$ of all antilinear maps $T : K \longrightarrow H$ such that $\sum_j \|Tu_j\|^2 < \infty$ for some, and hence every, orthonormal basis $E_2$ of $K$. Moreover for every $T \in H \otimes K$ we set*

$$\|\|T\|\|^2 = \sum_j \|Tu_j\|^2. \tag{5}$$

By Theorem 7.12 in [6], $H \otimes K$ is a Hilbert space with the norm $\|\|\cdot\|\|$ and associated inner product

$$<Q, T> = \sum_j <Qu_j, Tu_j>, \tag{6}$$

where $E_2 = \{u_j\}_{j \in J}$ is any orthonormal basis of $K$. Let $x \in H$ and $y \in K$. Then we define the map $x \otimes y$ by

$$(x \otimes y)(y') = <y, y'> x, \qquad (y' \in K). \tag{7}$$

Obviously $x \otimes y$ belongs to $H \otimes K$.

Let $T \in H \otimes K$. If $x, x' \in H$ and $y, y' \in K$, then by [6]

$$\|\|T\|\| = \|\|T^*\|\|, \tag{8}$$
$$\|\|x \otimes y\|\| = \|x\|\|y\|, \tag{9}$$
$$<x \otimes y, x' \otimes y'> = <x, x'><y, y'>. \tag{10}$$

Suppose $E_1 = \{e_i\}_{i \in I}$ and $E_2 = \{u_j\}_{j \in J}$ are orthonormal bases for $H$ and $K$, respectively. Then,

$$E_1 \otimes E_2 = \{e_i \otimes u_j : i \in I, j \in J\}$$

is an orthonormal basis for $H \otimes K$, by Proposition 7.14 in [6].

Now we can generalize Theorem 2.6.4 of [11] and Proposition 7.14 of [6] to frames.

**Theorem 2.2** *Let $H_1, \cdots, H_n$ be Hilbert spaces and $Y_1 = \{y_{1,i}\}_{i \in J_1}, \cdots, Y_n = \{y_{n,i}\}_{i \in J_n}$ be frames for $H_1, \cdots, H_n$, with frame bounds $A_1, B_1; \cdots; A_n, B_n$, respectively. Then*

$$\{y_{1,i_1} \otimes \cdots \otimes y_{n,i_n} \ : \ y_{j,i_j} \in Y_j, \ 1 \leq j \leq n\}$$

*is a frame for $H_1 \otimes \cdots \otimes H_n$ with frame bounds $A_1 A_2 \cdots A_n$ and $B_1 B_2 \cdots B_n$. In particular, if $Y_1, \cdots, Y_n$ are normalized tight frames, then it is a normalized tight frame.*

By using the associativity of tensor product [11, Proposition 2.6.5] and by induction it is enough to prove the theorem for $n = 2$.

**Theorem 2.3** *Let $\{x_n\}_{n \in I}$ and $\{y_m\}_{m \in J}$ be frames for $H$ and $K$, respectively. Then $\{x_n \otimes y_m\}_{n \in I, m \in J}$ is a frame for $H \otimes K$. Moreover, $\{x_n \otimes y_m\}$ is a normalized tight frame if $\{x_n\}$ and $\{y_m\}$ are.*

**Proof.** Let $A, B$ and $C, D$ be the bounds of the frames $\{x_n\}$ and $\{y_n\}$, respectively. Then by the Parseval identity, for all $T \in H \otimes K$ we have

$$<T, x_n \otimes y_m> = \sum_j <Tu_j, x_n \otimes y_m(u_j)> = \sum_j <Tu_j, <y_m, u_j> x_n>$$

$$= \sum_j \overline{<y_m, u_j>} <Tu_j, x_n>$$

$$= <\sum_j \overline{<y_m, u_j>} Tu_j, x_n> \qquad (T \text{ is an antilinear map})$$

$$= <T(\sum_j <y_m, u_j> u_j), x_n> = <Ty_m, x_n>.$$

Therefore $\sum_n \sum_m |<T, x_n \otimes y_m>|^2 = \sum_n \sum_m |<Ty_m, x_n>|^2$. Since $\{x_n\}$ is a frame for $H$, it follows that for every $m \in J_2$,

$$A\|Ty_m\|^2 \leq \sum_n |<Ty_m, x_n>|^2 \leq B\|Ty_m\|^2,$$

and

$$A \sum_m \|Ty_m\|^2 \leq \sum_m \sum_n |<T, x_n \otimes y_m>|^2 \leq B \sum_m \|Ty_m\|^2. \qquad (11)$$

Moreover, since $E_1$ is an orthonormal basis for $H$, then by the Parseval identity

$$\|Ty_m\|^2 = \sum_i |<Ty_m, e_i>|^2 = \sum_i |<T^*e_i, y_m>|^2. \qquad (12)$$

Now by using the fact that $\{y_m\}$ is a frame for $K$ and by (4), we conclude that

$$\sum_m \|Ty_m\|^2 = \sum_m \sum_i |<T^*e_i, y_m>|^2$$
$$= \sum_i \sum_m |<T^*e_i, y_m>|^2$$
$$\leq D \sum_i \|T^*e_i\|^2 = D \sum_j \|Tu_j\|^2 = D\|\|T\|\|^2,$$

and similarly,

$$\sum_m \|Ty_m\|^2 \geq C \sum_i \|T^*e_i\|^2 = C \sum_j \|Tu_j\|^2 = C\|\|T\|\|^2.$$

Thus

$$C\|\|T\|\|^2 \leq \sum_m \|Ty_m\|^2 \leq D\|\|T\|\|^2. \tag{13}$$

Therefore by inequalities (11) and (13), we get

$$AC\|\|T\|\|^2 \leq \sum_n \sum_m |<T, x_n \otimes y_m>|^2 \leq BD\|\|T\|\|^2.$$

Thus $\{x_n \otimes y_m\}$ is a frame for $H \otimes K$. ◇

For the converse we have the following result.

**Theorem 2.4** *Let $\{T_n\}_{n \in J}$ be a frame for $H \otimes K$. Then for each $x_0 \in H$ and $y_0 \in K$ the sequences $\{T_n y_0\}_{n \in J}$ and $\{T_n^* x_0\}_{n \in J}$ are frames for $H$ and $K$, respectively. Moreover these are tight frames, if $\{T_n\}_{n \in J}$ is.*

**Proof.** Let $A, B$ be the frame bounds for $\{T_n\}_{n \in J}$. As we saw in the proof of Theorem 2.1, for all $x \in H$ we have

$$<x \otimes y_0, T_n> = <x, T_n y_0>.$$

Since $\{T_n\}_{n \in J}$ is a frame for $H \otimes K$, we have

$$A\|x \otimes y_0\|^2 \leq \sum_n |<x \otimes y_0, T_n>|^2 \leq B\|x \otimes y_0\|^2.$$

Hence

$$A\|y_0\|^2\|x\|^2 \leq \sum_n |<x, T_n y_0>|^2 \leq B\|y_0\|^2\|x\|^2.$$

Therefore $\{T_n y_0\}_{n \in J}$ is a frame for $H$. Similarly, since for all $y \in K$

$$<y, T_n^* x_0> = <x_0, T_n y> = <x_0 \otimes y, T_n>,$$

we conclude that $\{T_n^* x_0\}_{n \in J}$ is also a frame for $K$. ◇

**Corollary 2.5** *If $\{T_n\}_{n \in J}$ is a frame for $H \otimes K$, then for each $x_0 \in H$ and $y_0 \in K$ the sequence $\{T_n(y_0 \otimes x_0)T_m\}$ is also a frame for $H \otimes K$.*

**Proof.** In view of Theorems 2.3 and 2.4, the sequence $\{T_n y_0 \otimes T_m^* x_0\}$ is a frame for $H \otimes K$. On the other hand, since $T_n$ is antilinear, for every $y \in K$

$$T_n y_0 \otimes T_m^* x_0 (y) = <T_m^* x_0, y> T_n y_0 = <T_m y, x_0> T_n y_0 = T_n(y_0 \otimes x_0) T_m(y).$$

Hence $\{T_n(y_0 \otimes x_0) T_m\}$ is a frame for $H \otimes K$. ◇

**Theorem 2.6** *If $Q \in B(H)$ is an invertible operator and $\{T_n\}_{n \in J}$ is a frame in $H \otimes K$, then the sequence $\{QT_n\}_{n \in J}$ is also a frame for $H \otimes K$.*

**Proof.** Since $Q$ is a bounded invertible operator on $H$, then for each $x \in H$

$$\|Q^{-1}\|^{-1} \|x\| \leq \|Q^* x\| \leq \|Q\| \|x\|. \tag{14}$$

Let $T \in H \otimes K$. Since $\{T_n\}$ is a frame for $H \otimes K$ and $Q^* T \in H \otimes K$ we have

$$A \|\|Q^* T\|\|^2 \leq \sum_n |<Q^* T, T_n>|^2 \leq B \|\|Q^* T\|\|^2,$$

where $A, B$ are frame bounds for $\{T_n\}_{n \in J}$. But $<Q^* T, T_n> = <T, QT_n>$, therefore

$$A \|\|Q^* T\|\|^2 \leq \sum_n |<T, QT_n>|^2 \leq B \|\|Q^* T\|\|^2. \tag{15}$$

Now by using (5) and (14) for every $j \in J$, we have

$$A \|Q^{-1}\|^{-2} \|\|T\|\|^2 \leq \sum_n |<T, QT_n>|^2 \leq B \|Q\|^2 \|\|T\|\|^2$$

Therefore $\{QT_n\}_{n \in J}$ is a frame for $H \otimes K$. ◇

**Corollary 2.7** *If $Q \in B(H)$ is a unitary operator and $\{T_n\}_{n \in J}$ is a frame in $H \otimes K$, then the sequence $\{R_n\}_{n \in J}$ defined by $R_n = QT_n$ is also a frame for $H \otimes K$.*

Since $\|\|T\|\| = \|\|T^*\|\|$, then $T \in H \otimes K$ if and only if $T^* \in K \otimes H$. Now we have the following result.

**Theorem 2.8** *The sequence $\{T_n\}_{n \in J}$ in the Hilbert space $H \otimes K$ is a frame if and only if $\{T_n^*\}_{n \in J}$ is a frame for $K \otimes H$.*

**Proof.** It's enough to note that $\|\|T^*\|\| = \|\|T\|\|$ and for evrey $n$, by applying the Parseval identity as in (4), we have $<T^*, T_n> = <T, T_n^*>$.

For the converse, it is enough to note that $T_n^{**} = T_n$. ◇

**Corollary 2.9** *The sequence $\{T_n\}_{n \in J}$ in the Hilbert space $H \otimes H$ is a frame if and only if $\{T_n^*\}_{n \in J}$ is.*

**Corollary 2.10** *If $Q \in B(K)$ is an invertible operator and $\{T_n\}_{n \in J}$ is a frame for $H \otimes K$, then the sequence $\{T_n Q\}_{n \in J}$ is also a frame for $H \otimes K$.*

**Proof.** Since $\{T_n\}_{n\in J}$ is a frame for $H \otimes K$, $\{T_n^*\}_{n\in J}$ is a frame for $K \otimes H$, by the above theorem. Moreover, $Q \in B(K)$ is an invertible operator, thus $Q^* \in B(K)$ is an invertible operator. Now by Theorem 2.6, $\{Q^*T_n^*\}_{n\in J}$ is a frame for $K \otimes H$. Therefore, by Theorem 2.8, $\{T_nQ\}_{n\in J}$ is also a frame for $H \otimes K$. ◇

**Corollary 2.11** *Let $\{T_n\}_{n\in J}$ be a frame for $H \otimes K$ and let $Q \in B(H)$ and $R \in B(K)$ be invertible operators. Then the sequence $\{S_n\}_{n\in J}$ defined by $S_n = QT_nR$ is a frame for $H \otimes K$.*

Now we can represent the inner product in $H \otimes K$ by normalized tight frames in $K$.

**Theorem 2.12** *If $\{y_m\}_{m\in J}$ is a normalized tight frame for Hilbert space $K$, then for all $Q, T \in H \otimes K$ we have*

$$<Q,T> = \sum_{m\in J} <Qy_m, Ty_m>.$$

**Proof.** Let $T, Q \in H \otimes K$. Let $\{e_i\}_{i\in I}$ be an orthonormal basis for $K$. Since $\{y_m\}$ is a normalized tight frame, then for every $i \in I$ and $T \in H \otimes K$, we have $T^*e_i = \sum_m <T^*e_i, y_m> y_m$ and $\|T^*e_i\|^2 = \sum_m |<T^*e_i, y_m>|^2$. Hence

$$<Q,T> = \sum_i <Q^*e_i, T^*e_i> = \sum_i <Q^*e_i, \sum_m <T^*e_i, y_m> y_m>$$
$$= \sum_i \sum_m <Q^*e_i, y_m> \overline{<Ty_m, e_i>}.$$

Now by using the Schwarz inequality, we conclude that the last double series is absolutely convergent, so

$$<Q,T> = \sum_m <Qy_m, \sum_i <Ty_m, e_i> e_i> = \sum_m <Qy_m, Ty_m>. \quad ◇$$

**Corollary 2.13** *If $T \in H \otimes K$ and $\{y_n\}_{n\in J_1}$, $\{z_m\}_{m\in J_2}$ be two normalized tight frames for $K$, then*

$$\sum_n \|Ty_n\|^2 = \sum_m \|Tz_m\|^2.$$

**Proposition 2.14** *Let $E_1 = \{e_i\}_{i\in I}$ and $E_2 = \{u_j\}_{j\in J}$ be orthonormal bases for $H$ and $K$, respectively, and let $T \in H \otimes K$. Then*

$$T = \sum_i e_i \otimes T^*e_i \quad \text{(pointwise)},$$
$$T = \sum_j Tu_j \otimes u_j \quad \text{(pointwise)}.$$

**Proof.** By the Parseval identity for every $y \in K$ we have

$$Ty = \sum_i <Ty, e_i> e_i = \sum_i <T^*e_i, y> e_i$$
$$= \sum_i e_i \otimes T^*e_i(y),$$

and

$$Ty = T(\sum_j <y, u_j> u_j) = \sum_j \overline{<y, u_j>} Tu_j$$
$$= \sum_j <u_j, y> Tu_j = \sum_j Tu_j \otimes u_j(y).$$

Thus $T = \sum_i e_i \otimes T^*e_i = \sum_j Tu_j \otimes u_j$.      ◇

**Proposition 2.15** *Let $\{x_n\}_{n \in J_1}$ and $\{y_m\}_{m \in J_2}$ be normalized tight frames for $H$ and $K$, respectively, and let $T \in H \otimes K$. Then $T = \sum_n x_n \otimes T^*x_n = \sum_m Ty_m \otimes y_m$.*

**Proof.** Let $y \in K$ be arbitrary. Then $y = \sum_m <y, y_m> y_m$ and $Ty = \sum_n <Ty, x_n> x_n$, in the weak sense. Therefore

$$Ty = \sum_n <Ty, x_n> x_n = \sum_n <T^*x_n, y> x_n = \sum_n x_n \otimes T^*x_n(y),$$

and

$$Ty = T(\sum_m <y, y_m> y_m) = \sum_m \overline{<y, y_m>} Ty_m$$
$$= \sum_m <y_m, y> Ty_m = \sum_m Ty_m \otimes y_m(y).$$

Thus $T = \sum_n x_n \otimes T^*x_n = \sum_m Ty_m \otimes y_m$.      ◇

## 3   The Canonical Dual Frame

Let $\{x_n\}_{n \in J}$ be a frame in the Hilbert space $H$. Then the operator $S : H \longrightarrow H$ defined by $Sx = \sum_n <x, x_n> x_n$, $(x \in H)$, is called the *frame operator* for $\{x_n\}$. By Theorem 2.1.3 in [9], $\{x_n\}_{n \in J}$ is a frame with frame bounds $A, B$ if and only if $S$ is a bounded linear operator with $AI \leq S \leq BI$ where $I$ denotes the identity operator on $H$. Moreover $S$ is a positive bounded linear invertible operator, and the sequence $\{S^{-1}x_n\}_{n\in J}$ is a frame with frame bounds $B^{-1}, A^{-1}$ for $H$. Every $x \in H$ can be written as

$$x = \sum_n <x, S^{-1}x_n> x_n = \sum_n <x, x_n> S^{-1}x_n.$$

Thus $S^{-1}$ is the frame operator of $\{S^{-1}x_n\}$ and

$$<x, S^{-1}x> = \sum_n |<S^{-1}x, x_n>|^2.$$

The frame defined by $x'_n = S^{-1}x_n$ $(n \in J)$ is called *dual frame* of $\{x_n\}_{n\in J}$ in the frame literature. Thus for every $x \in H$, we have

$$x = \sum_n <x, x'_n> x_n = \sum_n <x, x_n> x'_n. \tag{16}$$

Let $H$ and $K$ be Hilbert spaces. Then by Theorem 7.18 in [6] for all $Q, Q' \in B(H)$ and $T, T' \in B(K)$ we have
 (a) $Q \otimes T \in B(H \otimes K)$ and $\|Q \otimes T\| = \|Q\|\|T\|$
 (b) $(Q \otimes T)(x \otimes y) = Qx \otimes Ty$ for all $x \in H, y \in K$
 (c) $(Q \otimes T)(Q' \otimes T') = (QQ') \otimes (TT')$
 (d) If $Q \in B(H)$ and $T \in B(K)$ be invertible operators, then $Q \otimes T$ is an invertible operator and $(Q \otimes T)^{-1} = Q^{-1} \otimes T^{-1}$.

**Proposition 3.1** *Let $\{x_n\}_{n\in J_1}$ and $\{y_m\}_{m\in J_2}$ be frames in the Hilbert spaces $H$ and $K$ respectively, and let $S_1, S_2$ and $S$ be the frame operators of $\{x_n\}$, $\{y_m\}$ and $\{x_n \otimes y_m\}$, respectively. Then $S = S_1 \otimes S_2$.*

**Proof.** Let $T \in H \otimes K$ be arbitrary. Then $T = \sum_j Tu_j \otimes u_j$, by Proposition 2.14. Hence

$$S(T) = \sum_n \sum_m <T, x_n \otimes y_m> x_n \otimes y_m$$

$$= \sum_n \sum_m <\sum_j Tu_j \otimes u_j, x_n \otimes y_m> x_n \otimes y_m$$

$$= \sum_n \sum_m \sum_j (<Tu_j, x_n><u_j, y_m>) x_n \otimes y_m$$

$$= \sum_j (\sum_n <Tu_j, x_n> x_n) \otimes (\sum_m <u_j, y_m> y_m)$$

$$= \sum_j S_1(Tu_j) \otimes S_2(u_j) = \sum_j S_1 \otimes S_2(Tu_j \otimes u_j)$$

$$= S_1 \otimes S_2(\sum_j Tu_j \otimes u_j) = S_1 \otimes S_2(T).$$

Thus $S = S_1 \otimes S_2$. ◇

**Proposition 3.2** Let $\{x_n\}_{n \in I}$ and $\{y_m\}_{m \in J}$ be frames in the Hilbert spaces $H$ and $K$, respectively. If $0 \neq \lambda \in \mathbb{C}$, then

(i) $(x_n \otimes y_m)' = x_n' \otimes y_m'$

(ii) $(\lambda x_n)' = (\bar{\lambda})^{-1} x_n'$

(iii) $(x_n')' = x_n$

**Proof.** (i) Let $S_1, S_2$ and $S$ be the frame operators of $\{x_n\}$, $\{y_m\}$ and $\{x_n \otimes y_m\}$, respectively. Then by Proposition 3.1, $S = S_1 \otimes S_2$ and $S^{-1} = S_1^{-1} \otimes S_2^{-1}$. Thus

$$(x_n \otimes y_m)' = S^{-1}(x_n \otimes y_m) = (S_1^{-1} \otimes S_2^{-1})(x_n \otimes y_m)$$
$$= S_1^{-1} x_n \otimes S_2^{-1} y_m = x_n' \otimes y_m'.$$

(ii) Let $T$ be the frame operator for $\{\lambda x_n\}$. Then $T^{-1} = |\lambda|^{-2} S_1^{-1}$, where $S_1$ is the frame operator of $\{x_n\}$. Therefore

$$(\lambda x_n)' = T^{-1}(\lambda x_n) = |\lambda|^{-2} S_1^{-1}(\lambda x_n)$$
$$= |\lambda|^{-2} \lambda S_1^{-1}(x_n) = (\bar{\lambda})^{-1} x_n'.$$

(iii) Obvious.                                                                                                      ◇

**Proposition 3.3** Let $\{T_n\}_{n \in J}$ be a frame for $H \otimes K$. Then $(T_n^*)' = (T_n')^*$.

**Proof.** Let $S_1 \in B(H \otimes K)$ and $S_2 \in B(K \otimes H)$ be the frame operators of $\{T_n\}$ and $\{T_n^*\}$, respectively. Then for every $T \in H \otimes K$, we have

$$S_1(T) = \sum_n <T, T_n> T_n \quad \text{and} \quad S_2(T^*) = \sum_n <T^*, T_n^*> T_n^*.$$

Since for every $n \in J$, $<T^*, T_n^*> = <T, T_n>$ and since the map $T \longmapsto T^*$ is a linear operator, we conclude that

$$S_2(T^*) = \sum_n <T^*, T_n^*> T_n^* = \sum_n <T, T_n> T_n^*$$
$$= (\sum_n <T, T_n> T_n)^* = (S_1(T))^*.$$

So $S_2(T^*) = (S_1(T))^*$. Since $S_1^{-1}(T_n) = T_n'$, it follows that $T_n = S_1(T_n')$. By taking adjoints on both sides, we get $T_n^* = (S_1 n(T_n'))^* = S_2((T_n')^*)$. Hence

$$(T_n')^* = S_2^{-1}(T_n^*) = (T_n^*)'$$

Therefore, $(T_n')^* = (T_n^*)'$.                                                                                                   ◇

## Contents